\documentclass[twocolumn,amsthm]{autart}    % Enable this line and disable the 
\usepackage{graphicx}          % Include this line if your 
\usepackage{amsfonts}
\usepackage{amsmath}
\usepackage{amssymb}
\usepackage{booktabs}
\usepackage{tabularx}
\usepackage{array}
\usepackage{algorithm}
\usepackage{algpseudocode}
\usepackage{natbib}
\usepackage{xcolor}
\usepackage{url}
\usepackage[normalem]{ulem} 
\usepackage{array}
\usepackage{tabularx}
\newcolumntype{Y}{>{\centering\arraybackslash}X}

\newtheorem{theorem}{Theorem}

\newtheorem{corollary}{Corollary}

\theoremstyle{definition}
\newtheorem{definition}{Definition}

\theoremstyle{remark}
\newtheorem{remark}{Remark}

\definecolor{algcomment}{gray}{0.45}
\algrenewcommand\algorithmiccomment[1]{%
\hfill{\color{algcomment}\(\triangleright\) #1}}

\usepackage{amsmath}
\usepackage{array}
\usepackage{booktabs}
\usepackage{bm}

\newcommand{\phmat}[1]{%
  \begingroup
  \footnotesize
  \setlength{\arraycolsep}{1.2pt}%
  \renewcommand{\arraystretch}{1.25}%
  $\begin{bmatrix}
  #1
  \end{bmatrix}$%
  \endgroup
}

\begin{document}

\begin{frontmatter}

\title{PH-Symb : Sparse symbolic identification of nonlinear port-Hamiltonian
systems using Kolmogorov--Arnold networks%
\thanksref{footnoteinfo}}

\thanks[footnoteinfo]{This work has been supported by the EIPHI Graduate School
(contract ANR-17-EURE-0002) and by the Bourgogne-Franche-Comt\'e Region. Computations have been performed on the supercomputer facilities of the Mésocentre de calcul de Franche-Comté.}

\author[femto]{Achraf {El Messaoudi}},
\author[femto]{Karim Cherifi},
\author[femto]{Yann {Le Gorrec}},
\author[femto]{Yongxin Wu}

\address[femto]{Universit\'e Marie et Louis Pasteur, SUPMICROTECH, CNRS,
institut FEMTO-ST, F-25000 Besan\c{c}on, France}                            
          
                             % keyword list or with the 
                                          % help of the Automatica 
                                          % keyword wizard

\begin{abstract}  
The paper introduces PH-Symb, a novel symbolic structure-preserving identification method, that combines Kolmogorov Arnold networks and sparse regression to build explicit port-Hamiltonian (pH) models from data. Unlike other Machine Learning or neural pH approaches, PH-Symb recovers physically interpretable and explicit-symbolic components and enforces pH structure by construction. It also introduces a new normalization method that preserves the pH structure. It allows handling both identification from input-state-output (ISO) data and input-output (IO) data, and applications show it delivers interpretable models and lower errors than neural baselines on synthetic and real-world systems.
\end{abstract}

\end{frontmatter}

\section{Introduction}
\label{sec:intro}

Building mathematical models of nonlinear dynamical systems from measured data is a central problem in systems and control \citep{ljung1999}. Reliable models are required for simulation, analysis, state estimation, and model-based control. Machine-learning methods make it possible to identify nonlinear models from input-state or input-output data \citep{pillonetto2025}; however, the resulting generic black-box models lack physical interpretability. The port-Hamiltonian (pH) framework is of particular interest in this context. It explicitly describes stored energy, power-conserving interconnection, dissipation, and power-conjugate ports \citep{maschke1992,vanderschaft2014,vanderschaft2017}, and supports passivity-based control \citep{ortega2002}. Identifying a pH model requires respecting the underlying interconnection structure, the dissipation, the closure relation between state and co-state variables, and the passivity of the input-output mapping.

Existing structure-preserving approaches to pH learning can be classified according to whether the system state is measured. When state measurements are available, Gaussian-process (GPs) models incorporate Hamiltonian or passivity priors \citep{beckers2022}, while structured neural models learn pH dynamics from input-state data \citep{desai2021,eidnes2023, linares2026string} or input-state-output data \citep{cherifi2025}. When only inputs and outputs are measured, the state must be inferred from the data. This has been addressed using input-output GPs models \citep{ewering2025} and output-error pH models \citep{moradi2026} based on SUBNET \citep{beintema2023}.

However, preserving pH structure does not by itself provide interpretability. This is particularly limiting in the pH framework, where the Hamiltonian and the matrices encode distinct physical mechanisms. Neural pH models generally hide these components inside multilayer perceptrons (MLPs), whereas GP-based approaches describe them using kernel functions. In contrast, Kolmogorov-Arnold networks (KANs) make the learned representation easier to inspect by representing each edge function as a trainable spline \citep{liu2024kan}. PH-KAN extends KANs by embedding the pH structure into their architecture \citep{elmessaoudi2026phkan}, but inspecting the learned splines does not directly yield compact analytical models. Alternatively, Sparse Identification of Nonlinear Dynamical systems (SINDy) \citep{brunton2016} returns explicit equations. If the original governing equations are recovered, one important advantage is that the identified model remains valid for initial conditions outside the training domain. Recent work has extended SINDy to pH systems \citep{myhr2026siphy}. However, this method assumes a known constant interconnection matrix, as well as constant diagonal dissipation and input matrices, and needs state trajectories. A general method that constructs explicit symbolic expressions for all pH components, while considering both the case where the system states are measured and the case where only input-output data are available, has to be developed. The main contributions of this paper are the following :

\begin{itemize}
    \item We introduce PH-Symb, the first symbolic identification method that provides explicit expressions for all the pH components with extensions to handle both input-state-output (ISO) and input-output (IO) data.
    \item We show that PH-Symb can recover structure-preserving multivariate symbolic compositions from a library of univariate functions.
    \item We propose a normalization method for learning-based pH identification that preserves the pH structure and power balance and provides an exact transformation back to the original coordinates.
\end{itemize}

The remainder of the paper is organized as follows. Section~\ref{sec:prelim} formulates the identification problem for ISO and IO data. Section~\ref{sec:method} presents PH-Symb, including its symbolic construction, normalization, and training procedure. Section~\ref{sec:experiments} evaluates the method on synthetic and real-world data and compares it with state-of-the-art structure-preserving neural methods. Section~\ref{sec:conclusion} concludes the paper.

\emph{Notation.} 
For a matrix $A$, $A^\top$ denotes its transpose, while $A\succeq 0$ means that $A$ is positive semidefinite. The identity matrix of size $n$ is denoted by $I_n$. $\|\cdot\|_2$ and $\|\cdot\|_F$ denote the Euclidean and Frobenius norms, respectively. For a differentiable scalar function $H$, $\nabla H$ denotes its gradient.

\begin{figure*}[!t]
    \centering
    \includegraphics[width=\textwidth]{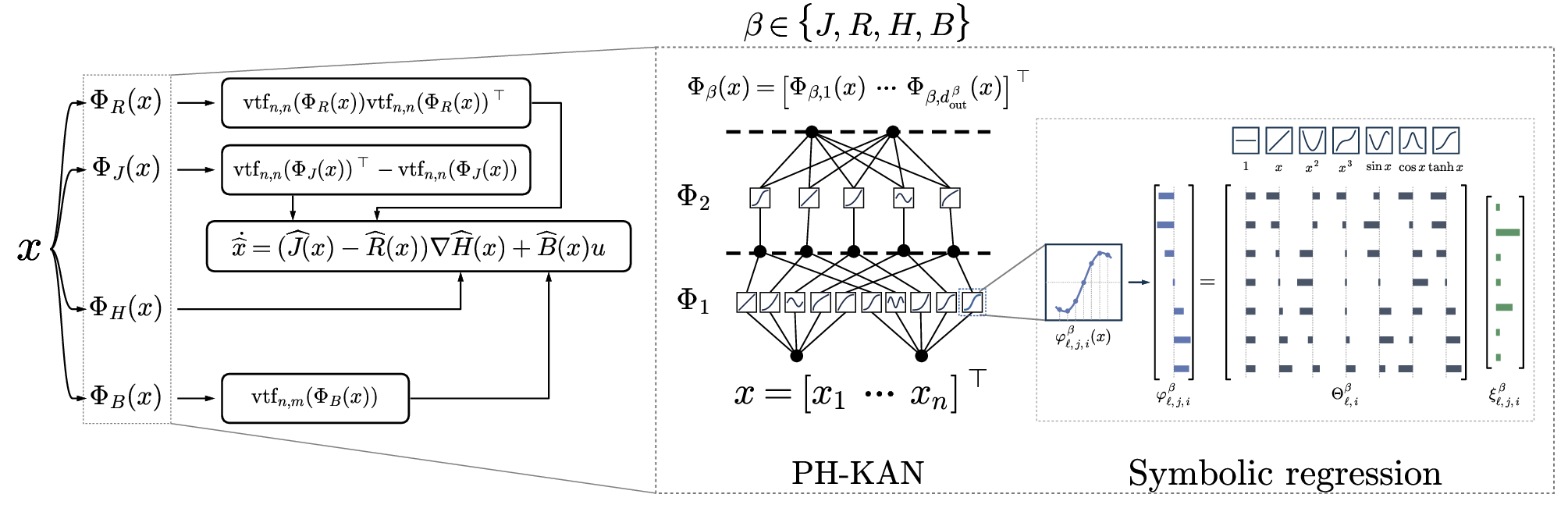}
    \caption{PH-Symb architecture combining PH-KAN with edge-wise symbolic regression. One PH-KAN block is assigned to each component $\Phi_\beta(x)$ where $\beta\in\{J,R,H,B\}$, and $\Phi_\beta(x)$ is then transformed and assembled into the structure-preserving pH model shown on the left. Each univariate KAN edge is represented using a prescribed symbolic library.}
    \label{fig:ph-kan-sindy}
\end{figure*}

\section{Problem formulation}
\label{sec:prelim}

\subsection{Port-Hamiltonian systems}
\label{sec:ph-systems}

We consider finite-dimensional nonlinear port-Hamiltonian systems of the form given in Definition~\ref{def:ph} \citep{maschke1992,vanderschaft2017}.

\begin{definition}[Nonlinear port-Hamiltonian system]
\label{def:ph}
Let $\mathcal X$ be a subset of $\mathbb R^n$. A port-Hamiltonian (pH) system with state $x(t)\in\mathcal X$, input $u(t)\in\mathbb R^m$, and output $y(t)\in\mathbb R^m$ is described by
\begin{equation}
\label{eq:ph}
\begin{aligned}
\dot x &= \bigl(J(x)-R(x)\bigr)\nabla H(x)+B(x)u,\\
y &= B(x)^\top\nabla H(x),
\end{aligned}
\end{equation}
where $H:\mathcal X\to\mathbb R$, a smooth positive definite function of $x$, is the Hamiltonian and $\nabla H$ denotes its gradient with respect to the state. The functions
\[
J,R:\mathcal X\to\mathbb R^{n\times n},
\qquad
B:\mathcal X\to\mathbb R^{n\times m}
\]
are the interconnection, dissipation, and input mapping matrices, respectively. For every $x\in\mathcal X$, $J$ and $R$ satisfy
\begin{equation}
\label{eq:constraints}
J(x)=-J(x)^\top,
\qquad
R(x)=R(x)^\top\succeq0.
\end{equation}
\end{definition}

The Hamiltonian describes energy storage. The matrix $B$ also defines the output dynamics in \eqref{eq:ph}. Consequently, the input and output are power-conjugate, and $y^\top u$ represents the power supplied to the system.

When no ambiguity arises, we omit the dependence on $x$ and write $J$, $R$, and $B$. Along solutions of \eqref{eq:ph}, the skew-symmetry of $J$ and the fact that $R$ is a positive semi-definite matrix gives
\begin{equation}
\label{eq:balance}
\frac{\mathrm d}{\mathrm dt}H\bigl(x(t)\bigr)
=
-\nabla H(x)^\top R(x)\nabla H(x)+y^\top u
\leq y^\top u.
\end{equation}
This balance states that the rate of change of stored energy equals the supplied power minus the dissipated power.

\subsection{Observation settings}
\label{sec:problem}
To identify the system, we collect data from $N_e$ experiments. For each experiment $r$, the state, input, and output samples are denoted by $x_k^{(r)}$, $u_k^{(r)}$, and $y_k^{(r)}$, respectively, with $k=0,\ldots,N_r$. Two measurement settings are considered, the input-state-output ISO setting and the IO setting.

\subsubsection{Input--state--output data.}
In the ISO setting, the data contain measurements of the input, state, and output:
\[
\mathcal D_{\mathrm{ISO}}
=
\left\{
\bigl(x_k^{(r)},u_k^{(r)},y_k^{(r)}\bigr)
\right\}_{r=1,\ldots,N_e;\,k=0,\ldots,N_r}.
\]
Assuming that the measured state variables satisfy a pH model of the form \eqref{eq:ph}, with measured input and output forming the associated power-conjugate pair, the objective is to identify symbolic expressions for the unknown components $(\widehat J,\widehat R, \widehat H, \widehat B)$ that define a pH model of the form \eqref{eq:ph}, satisfy \eqref{eq:constraints}, and predict the state and output trajectories.

\begin{remark}[]
\label{rem:output-scaling}
Several neural pH identification methods are trained using state data alone. In that case, the identified pH components are not unique. In particular, for a constant scalar rescaling, we obtain, for any $c>0$, the transformation
\[
J\longmapsto J/c,
\qquad
R\longmapsto R/c,
\qquad
H\longmapsto cH,
\]
which leaves the state dynamics unchanged. This scaling therefore does not affect the state prediction, but it becomes important when the identified components are intended to have a physical interpretation. We therefore include measurements of the conjugate output. Since the transformation above changes the output from $y$ to $cy$, a measured output can be used to determine $c$ and thus recover the true system components.
\end{remark}

\subsubsection{Input--output data.}
In the IO setting, the state is not measured, and the available data are
\[
\mathcal D_{\mathrm{IO}}
=
\left\{
\bigl(u_k^{(r)},y_k^{(r)}\bigr)
\right\}_{r=1,\ldots,N_e;\,k=0,\ldots,N_r}.
\]
The learned pH model is driven by the applied inputs and fitted to the measured outputs. Its state is denoted by $z\in\mathbb R^{n_z}$ to distinguish it from the measured state $x$ in the ISO setting. For a chosen model order $n_z$, the model is written as
\begin{equation}
\label{eq:io-ph}
\begin{aligned}
\dot z
&=
\bigl(\widehat J(z)-\widehat R(z)\bigr)
\nabla\widehat H(z)+\widehat B(z)u,\\
\widehat y
&=
\widehat B(z)^\top\nabla\widehat H(z),
\end{aligned}
\end{equation}
where
$\widehat H:\mathbb R^{n_z}\to\mathbb R$,
$\widehat J,\widehat R:\mathbb R^{n_z}\to\mathbb R^{n_z\times n_z}$,
and
$\widehat B:\mathbb R^{n_z}\to\mathbb R^{n_z\times m}$.
The pH constraints
$\widehat J(z)=-\widehat J(z)^\top$ and
$\widehat R(z)=\widehat R(z)^\top\succeq0$
are imposed, and the predicted output is computed directly from $\widehat B$ and $\widehat H$. Additional properties of $\widehat H$ can be imposed through its parameterization. For example, non-negativity can be imposed directly, while input-convex neural networks \citep{amos2017inputconvex} or input-convex KANs can enforce convexity \citep{deschatre2025inputconvex}. Such extensions are not considered here.

\section{The PH-Symb framework}
\label{sec:method}

Figure~\ref{fig:ph-kan-sindy} summarizes the PH-Symb framework. A separate KAN is used for each pH component $\Phi_\beta(x)$ where $\beta\in\{J,R,H,B\}$ and takes $x$ as input for ISO data or $z$ for IO data. The output vectors $\Phi_J$, $\Phi_R$, and $\Phi_B$ are passed to a reshaping operation, which arranges their entries into matrices used to construct $\widehat J$, $\widehat R$, and $\widehat B$. The scalar output $\Phi_H$ directly defines $\widehat H$. Combining these KANs with a reshaping operation and the transformations imposing skew-symmetry on $\widehat J$ and positive semidefiniteness on $\widehat R$ gives the PH-KAN architecture, as detailed in Section~\ref{sec:ph-kan-parameterization}. Symbolic regression is then applied to each spline edge function to obtain an explicit expression from the prescribed library.

\subsection{Matrix function approximation using KANs}
\label{sec:kan}
Multilayer perceptrons (MLPs) are widely used as universal function approximators. Kolmogorov--Arnold networks (KANs) provide an alternative universal approximation architecture with a different parameterization \citep{liu2024kan}. Both architectures consist of successive layers of nodes connected by edges. In an MLP, the edges carry scalar weights, while each node adds a bias and applies a fixed activation function. In a KAN, each edge carries a trainable univariate function, represented by a spline, and each node sums the outputs of its incoming edge functions. In Figure~\ref{fig:ph-kan-sindy}, the black dots denote the nodes, where these sums are performed, and the small boxes denote the edge functions.

For a fixed component $\beta$, the corresponding KAN block is a composition of the form
\begin{equation}
\label{eq:kan-block}
\Phi_\beta=\Phi_{L_\beta}\circ\cdots\circ\Phi_1,
\end{equation}
where $L_\beta$ is the number of layers. The block in Figure~\ref{fig:ph-kan-sindy} contains the two layers $\Phi_1$ and $\Phi_2$.

Let $\mathbf h=[h_1,\ldots,h_{d^\beta_{\ell-1}}]^\top$ be the $d^\beta_{\ell-1}$-node input of layer $\Phi_\ell$. The value returned by the $j$th node of $\Phi_\ell$ is
\begin{equation}
\label{eq:kan-layer}
\left[\Phi_\ell(\mathbf h)\right]_j
=
\sum_{i=1}^{d^\beta_{\ell-1}}
\phi^\beta_{\ell,j,i}(h_i),
\qquad
j=1,\ldots,d^\beta_\ell,
\end{equation}
where $\phi^\beta_{\ell,j,i}:\mathbb R\to\mathbb R$ is the trainable spline on the edge from input node $i$ to output node $j$. Thus, $\phi^\beta_{\ell,j,i}$ denotes one edge function, $\Phi_\ell$ one layer, and $\Phi_\beta$ the complete block. For an input state $x$, the block returns
\begin{equation}
\label{eq:kan-output}
\Phi_\beta(x)
=
\begin{bmatrix}
\Phi_{\beta,1}(x)&\cdots&\Phi_{\beta,d^\beta_{\mathrm{out}}}(x)
\end{bmatrix}^{\!\top},
\end{equation}
where $d^\beta_{\mathrm{out}}=d^\beta_{L_\beta}$ is the number of output nodes. 
\begin{remark}
    In the IO setting, $x$ is replaced by $z$.
\end{remark}
\subsection{PH-KAN parameterization}
\label{sec:ph-kan-parameterization}

PH-KAN assigns one KAN block $\Phi_\beta$ to each pH component $\beta\in\{J,R,H,B\}$ \citep{elmessaoudi2026phkan}. 
To construct the matrix-valued pH components from the KAN outputs, we use the vector-to-full-matrix operation introduced in \citep{cherifi2025}:
\begin{equation}
\label{eq:vtf}
\begin{aligned}
\operatorname{vtf}_{n,m}:\mathbb R^{nm}&\longrightarrow\mathbb R^{n\times m},\\
v=(v_i)_{i=1}^{nm}
&\longmapsto
\begin{bmatrix}
v_1 & v_{n+1} & \cdots & v_{(m-1)n+1}\\
v_2 & v_{n+2} & \cdots & v_{(m-1)n+2}\\
\vdots & \vdots & & \vdots\\
v_n & v_{2n} & \cdots & v_{mn}
\end{bmatrix}.
\end{aligned}
\end{equation}
The pH components are defined by
\begin{equation}
\label{eq:ph-assembly}
\begin{aligned}
\widehat J(x)
&=
\operatorname{vtf}_{n,n}\!\left(\Phi_J(x)\right)^\top
-
\operatorname{vtf}_{n,n}\!\left(\Phi_J(x)\right),\\
\widehat R(x)
&=
\operatorname{vtf}_{n,n}\!\left(\Phi_R(x)\right)
\operatorname{vtf}_{n,n}\!\left(\Phi_R(x)\right)^\top,\\
\widehat H(x)
&=\Phi_H(x),\\
\widehat B(x)
&=\operatorname{vtf}_{n,m}\!\left(\Phi_B(x)\right).
\end{aligned}
\end{equation}
This construction results in
$\widehat J(x)=-\widehat J(x)^\top$ and
$\widehat R(x)=\widehat R(x)^\top\succeq0$
for every $x$. The gradient $\nabla\widehat H(x)$ is obtained by automatic differentiation.

\begin{remark}[]
Available prior knowledge can be incorporated by fixing a known pH component or by restricting a block to be constant or to depend only on selected state variables.
\end{remark}

\subsection{Sparse symbolic expressions on the KAN edges}
\label{sec:sindykan}

PH-KAN provides a more interpretable representation than MLP-based pH models because its learned edge functions can be inspected individually. However, these functions are represented by splines and are not symbolic by construction. Standard post-hoc KAN extraction matches each trained spline to a single transformed function from a prescribed library. Since this matching is performed only after training, the learned spline is not constrained to align with any function in the library, and its replacement may alter the trained model. Moreover, using a single library function excludes sparse combinations of candidate functions \citep{bagrow2025}. Sparse identification of nonlinear dynamical systems (SINDy) \citep{brunton2016} produces such sparse combinations, but representing composed nonlinearities may require a large library containing the relevant compositions \citep{myhr2026siphy}.

Following SINDy-KAN \citep{howard2026}, PH-Symb applies sparse regression to each PH-KAN edge during training. To keep the computational and memory costs manageable, the regression is performed on a minibatch of $N_b$ states samples rather than on the complete data set. These are measured states in the ISO setting and model states generated along selected simulation windows in the IO setting.

Consider an edge $\phi^\beta_{\ell,j,i}$ and let $x\in\mathbb R$ denote
its scalar input. As illustrated in Figure~\ref{fig:ph-kan-sindy}, the
spline is approximated by a linear combination of candidate functions
$\theta_1,\ldots,\theta_{n_\Theta}$:
\begin{equation}
\label{eq:xi-edge}
\phi^{\beta,\Xi}_{\ell,j,i}(x)
=
\boldsymbol\theta(x)^\top\boldsymbol\xi^\beta_{\ell,j,i}
=
\sum_{q=1}^{n_\Theta}
\xi^\beta_{\ell,j,i,q}\theta_q(x),
\end{equation}
where
\[
\boldsymbol\theta(x)
=
\begin{bmatrix}
\theta_1(x)&\cdots&\theta_{n_\Theta}(x)
\end{bmatrix}^{\!\top}
\]
contains the candidate functions and
$\boldsymbol\xi^\beta_{\ell,j,i}\in\mathbb R^{n_\Theta}$ contains their
coefficients. For example, the library in Figure~\ref{fig:ph-kan-sindy} contains
$1$, $x$, $x^2$, $x^3$, $\sin(x)$, $\cos(x)$, and $\tanh(x)$.
The matrix $\boldsymbol\Xi$ collects the coefficient vectors
$\boldsymbol\xi^\beta_{\ell,j,i}$ row-wise over all KAN edges.

Over a minibatch, every edge receives $N_b$ scalar inputs $x(t_1),\ldots,x(t_{N_b})$. The candidate functions and the spline are evaluated at these inputs and collected in
\begin{equation}
\label{eq:edge-library}
\Theta^\beta_{\ell,i}
=
\begin{bmatrix}
\boldsymbol\theta(x(t_1))^\top\\
\vdots\\
\boldsymbol\theta(x(t_{N_b}))^\top
\end{bmatrix},
\qquad
\boldsymbol\varphi^\beta_{\ell,j,i}
=
\begin{bmatrix}
\phi^\beta_{\ell,j,i}(x(t_1))\\
\vdots\\
\phi^\beta_{\ell,j,i}(x(t_{N_b}))
\end{bmatrix}.
\end{equation}
Each row of $\Theta^\beta_{\ell,i}$ corresponds to one edge input and each column to one candidate function, while $\boldsymbol\varphi^\beta_{\ell,j,i}$ contains the corresponding spline values.

The coefficients are obtained by ridge regression:
\begin{equation}
\label{eq:edge-fit}
\boldsymbol\xi^\beta_{\ell,j,i}
=
\underset{\mathbf c\in\mathbb R^{n_\Theta}}
{\operatorname{arg\,min}}
\left\|
\boldsymbol\varphi^\beta_{\ell,j,i}
-
\Theta^\beta_{\ell,i}\mathbf c
\right\|_2^2
+
\gamma\|\mathbf c\|_2^2,
\qquad
\gamma>0,
\end{equation}
 where $\mathbf c\in\mathbb R^{n_\Theta}$ denotes any candidate vector of library coefficients, whereas $\boldsymbol\xi^\beta_{\ell,j,i}$ denotes the vector that minimizes the fitting error. The first term fits the symbolic expression to the spline values. The scalar $\gamma>0$ is the ridge regularization parameter. It penalizes large coefficients and stabilizes the solution when some columns of $\Theta^\beta_{\ell,i}$ are nearly dependent. This fit is recomputed as the splines and minibatch change during training.

The matrix $\boldsymbol\Xi$ obtained by ridge regression is generally
not sparse. A trainable matrix $\boldsymbol\Lambda$ of the same size is
therefore introduced to obtain sparse symbolic expressions. For each edge, the corresponding row $\boldsymbol\lambda^\beta_{\ell,j,i}\in\mathbb R^{n_\Theta}$ defines
\begin{equation}
\label{eq:symbolic-edges}
\phi^{\beta,\Lambda}_{\ell,j,i}(x)
=
\boldsymbol\theta(x)^\top
\boldsymbol\lambda^\beta_{\ell,j,i}
=
\sum_{q=1}^{n_\Theta}
\lambda^\beta_{\ell,j,i,q}\theta_q(x).
\end{equation}
The training of $\boldsymbol\Lambda$ and its coupling with $\boldsymbol\Xi$ are described in Section~\ref{sec:training}. After training, the retained edge expressions are composed through the KAN layers to construct the symbolic pH components.

\subsection{Structure-preserving normalization}
\label{sec:normalization}
Large differences in the numerical scales of the states, inputs and outputs can make training poorly conditioned, especially in multiphysics systems. In a pH model, these variables must be transformed together with the Hamiltonian, the system matrices, and the conjugate output to preserve the pH structure and power balance. Related transformations have been studied for linear and time-varying pH systems \citep{gugercin2012structure,cherifi2025dissipativity}. The following result gives the transformation used here for nonlinear pH systems.

\begin{theorem}[Structure-preserving normalization]
\label{thm:normalization}
Consider the pH system \eqref{eq:ph}. Let
$D\in\mathbb R^{n\times n}$ and
$S_u\in\mathbb R^{m\times m}$ be constant invertible matrices,
$x_0\in\mathbb R^n$, and $\tau,\kappa>0$. Define
\begin{equation}
\label{eq:norm-change}
\begin{aligned}
x&=x_0+D\widetilde x,
&\qquad
u&=S_u\widetilde u,\\
t&=\tau\widetilde t,
&\qquad
\widetilde H(\widetilde x)
&=\kappa H(x_0+D\widetilde x).
\end{aligned}
\end{equation}
For $x=x_0+D\widetilde x$, define
\begin{equation}
\label{eq:norm-maps}
\begin{aligned}
\widetilde J(\widetilde x)
&=
\frac{\tau}{\kappa}D^{-1}J(x)D^{-\top},\\
\widetilde R(\widetilde x)
&=
\frac{\tau}{\kappa}D^{-1}R(x)D^{-\top},\\
\widetilde B(\widetilde x)
&=
\tau D^{-1}B(x)S_u.
\end{aligned}
\end{equation}
Then the transformed system is port-Hamiltonian:
\begin{equation}
\label{eq:normalized-ph}
\begin{aligned}
\frac{\mathrm d\widetilde x}{\mathrm d\widetilde t}
&=
\bigl(
\widetilde J(\widetilde x)
-
\widetilde R(\widetilde x)
\bigr)
\nabla_{\widetilde x}\widetilde H(\widetilde x)
+
\widetilde B(\widetilde x)\widetilde u,\\
\widetilde y
&=
\widetilde B(\widetilde x)^\top
\nabla_{\widetilde x}\widetilde H(\widetilde x),
\end{aligned}
\end{equation}
where
$\widetilde J=-\widetilde J^\top$ and
$\widetilde R=\widetilde R^\top\succeq0$.
The transformation also preserves the power balance.
\end{theorem}

\begin{proof}
The chain rule gives
\[
\frac{\mathrm d\widetilde x}{\mathrm d\widetilde t}
=
\tau D^{-1}\dot x,
\qquad
\nabla_x H
=
\frac{1}{\kappa}
D^{-\top}\nabla_{\widetilde x}\widetilde H.
\]
Substitution into \eqref{eq:ph}, together with
$u=S_u\widetilde u$, directly gives
\eqref{eq:normalized-ph} with the maps in
\eqref{eq:norm-maps}.

Since $\tau/\kappa>0$, the congruence transformations in
\eqref{eq:norm-maps} preserve the skew-symmetry of $J$ and the
symmetry and positive semidefiniteness of $R$. Consequently,
\[
\frac{\mathrm d\widetilde H}{\mathrm d\widetilde t}
=
-
\nabla_{\widetilde x}\widetilde H^\top
\widetilde R
\nabla_{\widetilde x}\widetilde H
+
\widetilde y^\top\widetilde u.
\]
Finally,
\[
\frac{\mathrm d\widetilde H}{\mathrm d\widetilde t}
=
\tau\kappa\frac{\mathrm dH}{\mathrm dt},
\qquad
\widetilde y^\top\widetilde u
=
\tau\kappa y^\top u,
\]
so the complete power balance is scaled by the same positive factor
$\tau\kappa$.
\end{proof}

\begin{corollary}[Inverse transformation]
\label{cor:pipeline}
Let $\widetilde x=D^{-1}(x-x_0)$. The pH components and output in the
original coordinates are recovered as
\begin{equation}
\label{eq:inverse-normmaps}
\begin{aligned}
H(x)
&=
\frac{1}{\kappa}\widetilde H(\widetilde x),
&
J(x)
&=
\frac{\kappa}{\tau}
D\widetilde J(\widetilde x)D^\top,
\\
R(x)
&=
\frac{\kappa}{\tau}
D\widetilde R(\widetilde x)D^\top,
&
B(x)
&=
\frac{1}{\tau}
D\widetilde B(\widetilde x)S_u^{-1},
\\
y
&=
\frac{1}{\tau\kappa}
S_u^{-\top}\widetilde y.
\end{aligned}
\end{equation}
\end{corollary}

Thus, symbolic expressions identified in normalized coordinates can be
mapped exactly back to the original coordinates.

Let $\operatorname{std}_{\mathrm{tr}}(v)$ denote the standard deviation
of $v$ over the training data. In the ISO experiments, we set
$x_0=0$, $\kappa=1$, and
\[
\begin{aligned}
D
&=
\operatorname{diag}\!\left(
\operatorname{std}_{\mathrm{tr}}(x_1),\ldots,
\operatorname{std}_{\mathrm{tr}}(x_n)
\right),
\\
S_u
&=
\operatorname{diag}\!\left(
\operatorname{std}_{\mathrm{tr}}(u_1),\ldots,
\operatorname{std}_{\mathrm{tr}}(u_m)
\right).
\end{aligned}
\]
For state-derivative targets $\zeta_1,\ldots,\zeta_N$, the time scale is
\begin{equation}
\label{eq:empirical-time-scale}
\tau
=
\left[
\operatorname{std}_{\mathrm{tr}}\!\left(
\operatorname{vec}\!\left(
\{D^{-1}\zeta_k\}_{k=1}^{N}
\right)
\right)
\right]^{-1}.
\end{equation}
Here, $\operatorname{vec}$ stacks the normalized derivative vectors.
All standard deviations are computed from the training data and
assumed finite and strictly positive.

In the IO experiments, the measured signals are centered and
standardized as
\begin{equation}
\label{eq:deviation-ports}
\widetilde u=S_u^{-1}(u-u_0),
\qquad
\widetilde y=S_y^{-1}(y-y_0).
\end{equation}
These signals define the ports of the identified IO model. Since its
internal state is learned and $u_0$ and $y_0$ are external centering
constants, \eqref{eq:deviation-ports} is not interpreted as a coordinate
transformation of an unknown physical pH realization. Without offsets,
compatibility with Theorem~\ref{thm:normalization} requires
\begin{equation}
\label{eq:dual-port-scaling}
S_y^{-1}
=
\tau\kappa S_u^\top.
\end{equation}
In the SISO case, $\kappa$ can be chosen to satisfy this relation for
any positive input and output scales. 

\begin{remark}
Since Theorem~\ref{thm:normalization} follows directly from the pH state and output equations in \eqref{eq:ph}, the same normalization can be used with other learning-based pH identification methods.
\end{remark}

\begin{table*}[!t]
\centering
\caption{Ground-truth and PH-Symb models for the three benchmarks.}
\label{tab:identified-ph-components}
\scriptsize
\setlength{\tabcolsep}{0.8pt}
\renewcommand{\arraystretch}{1.35}
\begin{tabular*}{\textwidth}{
@{\extracolsep{\fill}}
c cc cc cc
@{}}
\toprule
&
\multicolumn{2}{c}{\textbf{Duffing}}
&
\multicolumn{2}{c}{\textbf{PMSM}}
&
\multicolumn{2}{c}{\textbf{MagLev}}
\\
\cmidrule(lr){2-3}
\cmidrule(lr){4-5}
\cmidrule(lr){6-7}
&
\textbf{True} & \textbf{PH-Symb}
&
\textbf{True} & \textbf{PH-Symb}
&
\textbf{True} & \textbf{PH-Symb}
\\
\midrule
$\bm{H(x)}$
&
$\begin{gathered}
0.5q^2\\
{}+0.25q^4\\
{}+0.5p^2
\end{gathered}$
&
$\begin{gathered}
0.4742q^2\\
{}+0.2393q^4\\
{}+0.4980p^2
\end{gathered}$
&
$\begin{gathered}
131.58\phi_d^2\\
{}+131.58\phi_q^2\\
{}+41.67p^2
\end{gathered}$
&
$\begin{gathered}
136.43\phi_d^2\\
{}+127.35\phi_q^2\\
{}+25.91p^2
\end{gathered}$
&
$\begin{gathered}
5p^2\\
{}+0.5q^2\phi^2\\
{}+0.05\phi^2
\end{gathered}$
&
$\begin{gathered}
4.877p^2\\
{}+0.510q^2\phi^2\\
{}+0.049\phi^2
\end{gathered}$
\\
\midrule
$\bm{J(x)}$
&
\phmat{
0&1\\
-1&0
}
&
\phmat{
0&1.0\\
-1.0&0
}
&
\phmat{
0&0&\phi_q\\
0&0&-(\phi_d+0.17)\\
-\phi_q&\phi_d+0.17&0
}
&
\phmat{
0&0.30p&0.05\phi_q\\
-0.30p&0&-0.21\\
-0.05\phi_q&0.21&0
}
&
\phmat{
0&1&0\\
-1&0&0\\
0&0&0
}
&
\phmat{
0&0.970&0\\
-0.970&0&0\\
0&0&0
}
\\
\addlinespace[5pt]
$\bm{R(x)}$
&
\phmat{
0&0\\
0&0.3
}
&
\phmat{
0&0\\
0&0.29
}
&
\phmat{
0.225&0&0\\
0&0.225&0\\
0&0&0.0026
}
&
\phmat{
0.217&0&0\\
0&0.235&0.052\\
0&0.052&0.0117
}
&
\phmat{
0&0&0\\
0&0.1&0\\
0&0&10
}
&
\phmat{
0&0&0\\
0&0.103&0\\
0&0&9.6
}
\\
\addlinespace[5pt]
$\bm{B(x)}$
&
\phmat{
0\\
1
}
&
\phmat{
0\\
1.00
}
&
\phmat{
1&0\\
0&1\\
0&0
}
&
\phmat{
1.00&0\\
0&1.00\\
0&0
}
&
\phmat{
0\\
0\\
1
}
&
\phmat{
0\\
0\\
0.985
}
\\
\addlinespace[5pt]
$\bm{y(x)}$
&
$p$
&
$0.9960p$
&
\phmat{
263.16\phi_d\\
263.16\phi_q
}
&
\phmat{
272.9\phi_d\\
254.7\phi_q
}
&
$\begin{gathered}
q^2\phi\\
{}+0.1\phi
\end{gathered}$
&
$\begin{gathered}
1.004q^2\phi\\
{}+0.096\phi
\end{gathered}$
\\
\bottomrule
\end{tabular*}
\end{table*}

\subsection{Training with ISO and IO data}
\label{sec:training}

\begin{figure}[!t]
    \centering
    \includegraphics[width=\columnwidth]{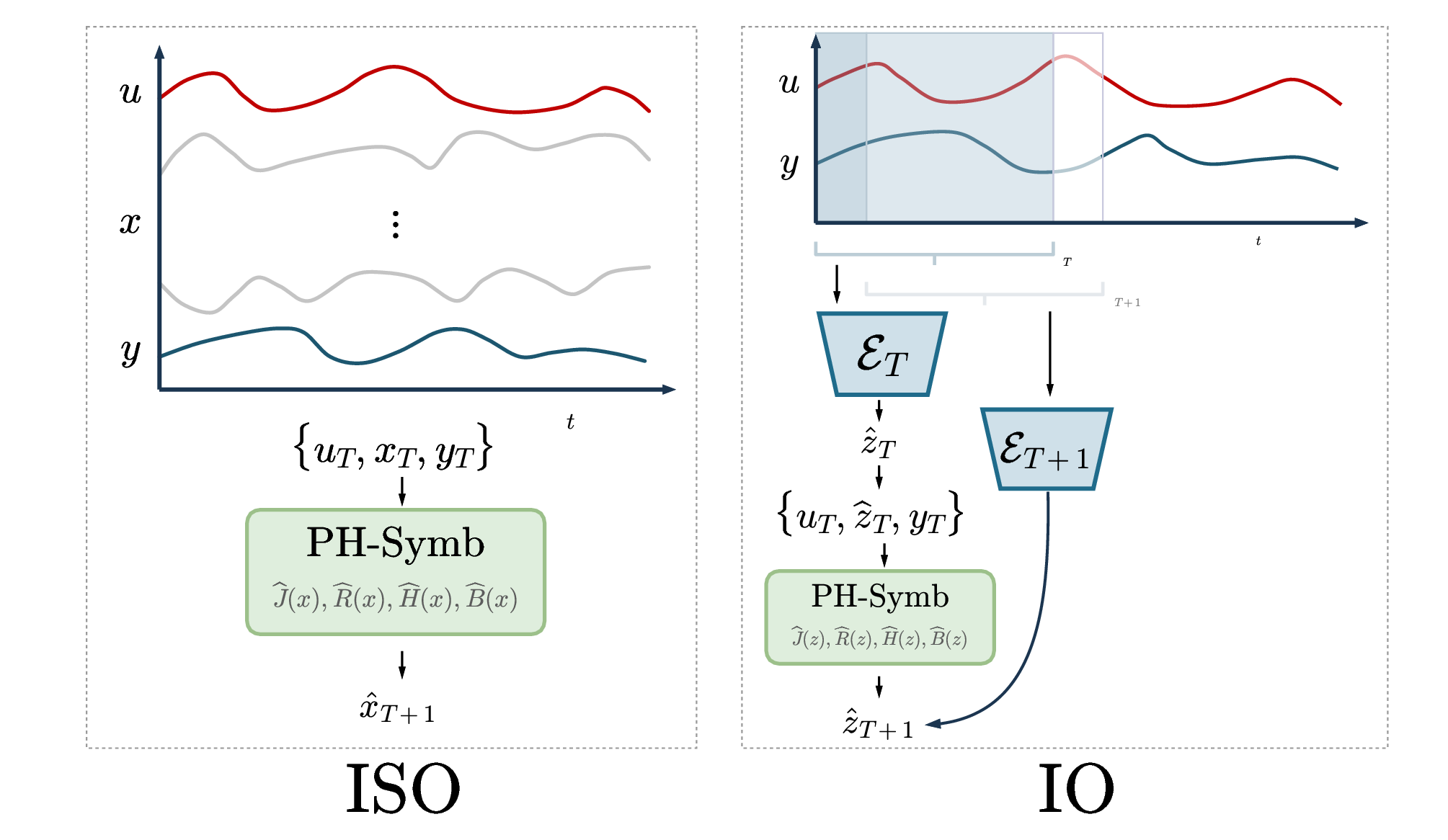}
    \caption{PH-Symb training with ISO and IO data. In the IO setting, the same encoder is evaluated on the initial and shifted input--output histories.}
    \label{fig:architecture}
\end{figure}

Figure~\ref{fig:architecture} summarizes the two training settings considered in this paper: ISO and IO.
The losses below are written using the normalized or standardized variables, with the tildes omitted for readability. During training, the data loss is evaluated for the spline model and the two symbolic models (dense and sparse), indexed by
$\nu\in\{\mathrm{KAN},\Xi,\Lambda\}$ respectively. The two symbolic representations are coupled through
\begin{equation}
\label{eq:symbolic-coupling}
\mathcal R_{\mathrm{cpl}}
=
\alpha_2
\left\lVert
\boldsymbol\Lambda-\boldsymbol\Xi
\right\rVert_F^2,
\qquad
\alpha_2\geq0.
\end{equation}
This term keeps the trainable symbolic coefficients $\boldsymbol\Lambda$ close to the coefficients $\boldsymbol\Xi$ obtained by fitting the current KAN edges. The data losses for the ISO and IO settings are defined next.

\subsubsection{ISO data}
\label{sec:iso-training}

In the ISO setting, the measured states $x_k$ are supplied directly to
the pH model. The derivative targets
$\zeta_k\simeq\dot x_k$ are estimated along each measured trajectory
using a Savitzky--Golay filter \citep{savitzky1964}. For a minibatch
$\mathcal B$ and representation $\nu\in\{\mathrm{KAN},\Xi,\Lambda\}$, define
\begin{equation}
\label{eq:iso-source-loss}
\begin{aligned}
\mathcal L^{\mathrm{ISO}}_\nu
={}&
\frac{\lambda_x}{n|\mathcal B|}
\sum_{k\in\mathcal B}
\left\|
\widehat f_\nu(x_k,u_k)-\zeta_k
\right\|_2^2
\\
&+
\frac{\lambda_y}{m|\mathcal B|}
\sum_{k\in\mathcal B}
\left\|
\widehat y_\nu(x_k)-y_k
\right\|_2^2,
\end{aligned}
\end{equation}
where $\widehat f_\nu$ and $\widehat y_\nu$ are the pH dynamics and
output obtained with representation
$\nu$, and
$\lambda_x,\lambda_y\geq0$. The ISO resulting objective is
\begin{equation}
\label{eq:iso-loss}
\mathcal L^{\mathrm{ISO}}
=
\sum_{\nu\in\{\mathrm{KAN},\Xi,\Lambda\}}
\omega_\nu\mathcal L^{\mathrm{ISO}}_\nu
+
\mathcal R_{\mathrm{cpl}}
+
\alpha_1\|\boldsymbol\Lambda\|_1,
\end{equation}
where $\omega_\nu\geq0$ weights the three representations and
$\alpha_1\geq0$ promotes sparsity. Output measurements may be included
during training by taking $\lambda_y>0$, or used after training to
adjust the constant scaling described in
Remark~\ref{rem:output-scaling}.

\subsubsection{IO data}
\label{sec:io-training}

In the IO setting, the state is not measured and the initial condition of each simulation window is unknown. As in SUBNET \citep{beintema2023} and the pH output-error approach \citep{moradi2026}, we use an encoder 
$\mathcal E$ to initialize each window from past input and output samples. Beyond this initialization that was introduced in \citep{beintema2023}, we apply the same encoder $\mathcal E$ at shifted positions within the window and introduces a state-consistency loss between the encoded and propagated states.

For a window starting at sample $k$ of experiment $r$, define the past output and input histories as
\begin{equation}
\label{eq:io-encoder}
\begin{aligned}
\mathbf y_{k,-}^{(r)}
&=
\begin{bmatrix}
\bigl(y_{k-n_a}^{(r)}\bigr)^\top &
\cdots &
\bigl(y_{k-1}^{(r)}\bigr)^\top
\end{bmatrix}^{\top},
\\
\mathbf u_{k,-}^{(r)}
&=
\begin{bmatrix}
\bigl(u_{k-n_b}^{(r)}\bigr)^\top &
\cdots &
\bigl(u_{k-1}^{(r)}\bigr)^\top
\end{bmatrix}^{\top},
\\
z_{\mathrm{enc},k}^{(r)}
&=
\mathcal{E}^{\eta}\!\left(
\mathbf y_{k,-}^{(r)},\mathbf u_{k,-}^{(r)}
\right).
\end{aligned}
\end{equation}
Here, $n_a$ and $n_b$ are the output- and input-history lengths. The encoder $\mathcal{E}^{\eta}$ is a feedforward neural network with parameters $\eta$ and $n_z$ outputs, which provide the initial model state $z_{\mathrm{enc},k}^{(r)}$ of length $n_z$.

Let $\mathcal B$ be a minibatch of pairs $(r,k)$ identifying an experiment and the beginning of a simulation window, and let $N_p$ be the prediction horizon. For each representation $\nu\in\{\mathrm{KAN},\Xi,\Lambda\}$, the model is initialized at $z_{\mathrm{enc},k}^{(r)}$ and simulated under the measured input. The state is propagated using the classical fourth-order Runge--Kutta method, with one integration step per sample and a zero-order hold on the input. We denote by $\widehat y_{k+j\mid k}^{(r),\nu}$ the output at sample $k+j$ predicted from a window initialized at sample $k$. Thus, $j$ is the position within the window and $k+j$ is the corresponding sample index. The output-prediction loss is
\begin{equation}
\label{eq:rollout-source-loss}
\mathcal L_\nu^{\mathrm{IO}}
=
\frac{1}{mN_p|\mathcal B|}
\sum_{(r,k)\in\mathcal B}
\sum_{j=0}^{N_p-1}
\left\|
\widehat y_{k+j\mid k}^{(r),\nu}
-
y_{k+j}^{(r)}
\right\|_2^2.
\end{equation}

To enforce state consistency, the encoder is evaluated again using the histories preceding sample $k+j$. Let $\widehat z_{k+j\mid k}^{(r),\Lambda}$ denote the state reached by the $\Lambda$ model after $j$ simulation steps from $z_{\mathrm{enc},k}^{(r)}$. For a set of offsets $\mathcal I\subseteq\{1,\ldots,N_p-1\}$, we define
\begin{equation}
\label{eq:encoder-alignment}
\mathcal L_{\mathrm{align}}
=
\frac{1}{n_z|\mathcal I||\mathcal B|}
\sum_{(r,k)\in\mathcal B}
\sum_{j\in\mathcal I}
\left\|
z_{\mathrm{enc},k+j}^{(r)}
-
\widehat z_{k+j\mid k}^{(r),\Lambda}
\right\|_2^2.
\end{equation}
The two terms inside the norm estimate the same internal state at sample $k+j$. The first is obtained directly from the shifted input--output history, while the second is propagated from the beginning of the window by the pH model.

The IO objective is
\begin{equation}
\label{eq:io-loss}
\mathcal L^{\mathrm{IO}}
=
\sum_{\nu\in\{\mathrm{KAN},\Xi,\Lambda\}}
\omega_\nu\mathcal L_\nu^{\mathrm{IO}}
+
\mathcal R_{\mathrm{cpl}}
+
\alpha_{\mathrm{align}}\mathcal L_{\mathrm{align}},
\end{equation}
where $\alpha_{\mathrm{align}}\geq0$ weights the state-consistency loss. Additional regularization terms, when used, are specified for each experiment.

After training, symbolic terms are selected according to their contribution over the inputs encountered by the corresponding edge. For the term $\lambda^\beta_{\ell,j,i,q}\theta_q(\zeta)$, define
\begin{equation}
\label{eq:coefficient-contribution}
a^\beta_{\ell,j,i,q}
=
\left|
\lambda^\beta_{\ell,j,i,q}
\right|
\left(
\frac{1}{N_{\mathrm{obs}}}
\sum_{p=1}^{N_{\mathrm{obs}}}
\theta_q(\zeta_p)^2
\right)^{1/2},
\end{equation}
where $\zeta_1,\ldots,\zeta_{N_{\mathrm{obs}}}$ are the values of the corresponding edge input encountered during training. Terms satisfying $a^\beta_{\ell,j,i,q}<\varepsilon$ are removed before assembling the symbolic pH components.

The complete training procedure is summarized in Algorithm~\ref{alg:training}.

\begin{algorithm}[!t]
\caption{PH-Symb training}
\label{alg:training}
\footnotesize
\begin{algorithmic}[1]
\Statex \textbf{Input:} Normalized ISO data or standardized IO data; derivative targets $\zeta$ for ISO; symbolic library; edge-regression minibatch size $N_b$; number of iterations $E$; threshold $\varepsilon$.
\State Initialize the four PH-KAN blocks and $\boldsymbol\Lambda$; for IO data, also initialize $\mathcal{E}^{\eta}$.
\For{$e=1,\ldots,E$}
\State Draw a minibatch $\mathcal B$.
\State Form a state minibatch $\mathcal S_{\mathcal B}$ of $N_b$ samples from $\{x_k^{(r)}\}_{(r,k)\in\mathcal B}$ for ISO data, or from $\{\widehat z_{k+j\mid k}^{(r),\Lambda}\}_{(r,k)\in\mathcal B,\,j=0,\ldots,N_p-1}$ obtained using \eqref{eq:io-encoder} and \eqref{eq:io-ph} for IO data.
\For{each KAN edge $\phi^\beta_{\ell,j,i}$}
\State Construct $\Theta^\beta_{\ell,i}$ and $\boldsymbol\varphi^\beta_{\ell,j,i}$ using \eqref{eq:edge-library}.
\State Fit the edge: compute $\boldsymbol\xi^\beta_{\ell,j,i}$ by solving \eqref{eq:edge-fit}.
\EndFor
\State Assemble and evaluate the KAN, $\boldsymbol\Xi$, and $\boldsymbol\Lambda$ pH models using \eqref{eq:ph-assembly}.
\State Set $\mathcal L\gets\mathcal L^{\mathrm{ISO}}$ using \eqref{eq:iso-loss}, or $\mathcal L\gets\mathcal L^{\mathrm{IO}}$ using \eqref{eq:io-loss}.
\State Gradient step: update the PH-KAN parameters and $\boldsymbol\Lambda$ using $\nabla\mathcal L$; for IO data, also update $\eta$.
\EndFor
\State Threshold $\boldsymbol\Lambda$ at $\varepsilon$ using $|\lambda^\beta_{\ell,j,i,q}|$ for ISO data and $a^\beta_{\ell,j,i,q}$ from \eqref{eq:coefficient-contribution} for IO data.
\State Assemble the final symbolic components $(\widehat J,\widehat R,\widehat H,\widehat B)$ from the thresholded $\boldsymbol\Lambda$ using \eqref{eq:ph-assembly}.
\State For ISO data, apply the inverse transformation \eqref{eq:inverse-normmaps}.
\State \Return $(\widehat J,\widehat R,\widehat H,\widehat B)$ and, for IO data, $\mathcal{E}^{\eta}$.
\end{algorithmic}
\end{algorithm}

\section{Numerical experiments}
\label{sec:experiments}

We now evaluate PH-Symb on ISO and IO benchmarks. In all experiments, the only prior knowledge is which pH components are linear or nonlinear and which nonlinear function types, such as polynomial, exponential, or square-root terms, should be included in the candidate library. No component is fixed in advance. The expressions and coefficients of $J$, $R$, $H$, and $B$ are learned in both the ISO and IO settings.

\begin{figure*}[!htpb]
  \centering
  \includegraphics[width=0.8\textwidth]{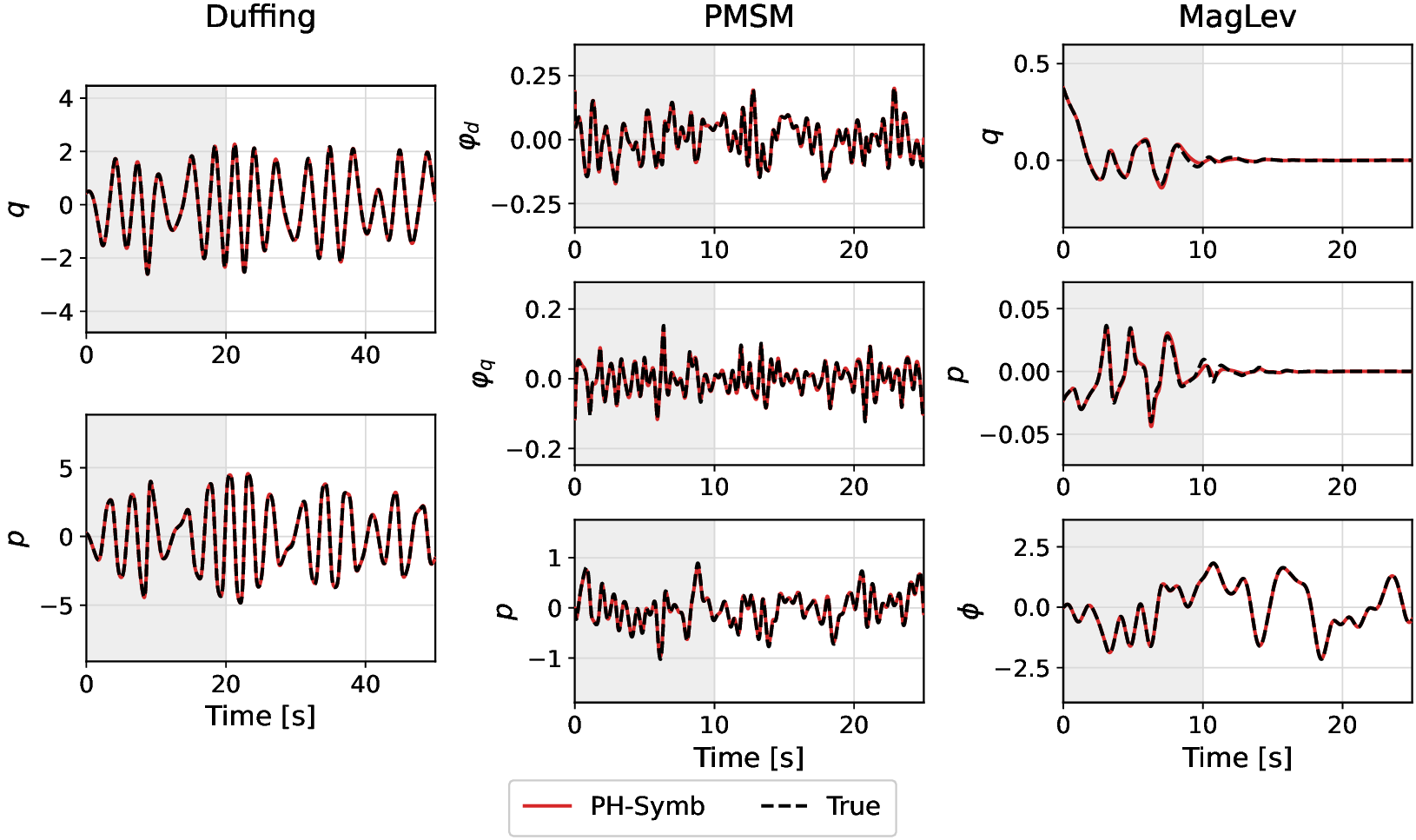}
  \caption{Representative free-run state predictions for the three ISO benchmarks. The shaded interval has the same duration as the training trajectories.}
  \label{fig:iso-rollouts}
\end{figure*}

\subsection{Input--state--output benchmarks}
\label{sec:iso-benchmarks}

We consider three benchmarks with different time scales and nonlinear port-Hamiltonian structures. The Duffing oscillator is a relatively slow mechanical system with a quartic Hamiltonian \citep{linares2026}. Its state is \(x=(q,p)\), where \(q\) and \(p\) are the position and momentum. The permanent-magnet synchronous motor (PMSM) combines faster electrical and mechanical dynamics and has a state-dependent interconnection matrix \citep{spirito2024}. Its state is \(x=(\phi_d,\phi_q,p)\), where \(\phi_d\) and \(\phi_q\) are the stator fluxes in the rotating \(dq\) frame and \(p\) is the rotor angular momentum. The magnetic levitation (MagLev) system also combines electrical and mechanical dynamics, but has a nonseparable Hamiltonian containing the multivariate product \(q^2\phi^2\) \citep{beckers2022}. Its state is \(x=(q,p,\phi)\), where \(q\) is the vertical position, \(p\) is the mechanical momentum, and \(\phi\) is the coil flux. The MagLev benchmark assesses the ability of PH-Symb to recover multivariate composite expressions from univariate edge libraries. For each benchmark, we use the physical model parameter values reported in the corresponding reference. The true pH components for each system can be found in Table~\ref{tab:identified-ph-components}.

For each input channel $\ell\in\{1,\ldots,m\}$, we use a random-phase multisine excitation \citep{pintelon2012},
\begin{equation}
u_\ell(t)
=
u_{\mathrm{RMS}}
\sqrt{\frac{2}{N_h}}
\sum_{h=1}^{N_h}
\sin\!\left(
\frac{2\pi h t}{T}
+
\varphi_{\ell h}
\right).
\label{eq:iso-input}
\end{equation}
Here, $N_h$ is the number of excited harmonics and $\varphi_{\ell h}$ is the phase of harmonic $h$ in input channel $\ell$, drawn independently from the uniform distribution on $[0,2\pi)$. The factor $u_{\mathrm{RMS}}\sqrt{2/N_h}$ sets the RMS value of the input to $u_{\mathrm{RMS}}$. The record duration $T$ determines the frequency grid, whose fundamental frequency is $1/T$. Duffing uses one input channel with $u_{\mathrm{RMS}}=1.6$ and $N_h=20$. PMSM uses two independently generated input channels with $u_{\mathrm{RMS}}=4.5$ and $N_h=40$. MagLev uses one input channel with $u_{\mathrm{RMS}}=2.0$ and $N_h=20$. The training and validation records last $20\,\mathrm{s}$ for Duffing and $10\,\mathrm{s}$ for PMSM and MagLev. The corresponding test records last $50\,\mathrm{s}$ and $25\,\mathrm{s}$, respectively. All signals are sampled with a period of $0.01\,\mathrm{s}$.

The initial states are drawn independently from uniform distributions. For Duffing, $q,p\in[-0.6,0.6]$. For PMSM, $\phi_d,\phi_q\in[-0.5,0.5]$ and $p\in[-1,1]$. For MagLev, $q,\phi\in[-1,1]$ and $p\in[-0.1,0.1]$. State derivatives used during training are estimated along the measured trajectories with a Savitzky--Golay filter.

Simulation accuracy is evaluated over the complete test record by initializing the identified model at the true initial state and applying the test input. For test experiment $r$ and state component $j$, the normalized root-mean-square error (NRMSE) is
\begin{equation}
e_{r,j}
=
\frac{
\left[
N_r^{-1}
\sum_{k=0}^{N_r-1}
\left(
x_{k,j}^{(r)}
-
\hat{x}_{k,j}^{(r)}
\right)^2
\right]^{1/2}
}{
\operatorname{std}_{\mathcal{D}_{\mathrm{train}}}(x_j)
}.
\label{eq:iso-nrmse}
\end{equation}
The normalization factor is computed from the clean training data. The reported errors are $100e_{r,j}$, expressed as percentages, and summarized by their mean and standard deviation over ten independent test experiments.

\subsection{ISO results}
\label{sec:iso-results}

Table~\ref{tab:identified-ph-components} compares the symbolic pH components recovered by PH-Symb with the ground-truth components. For Duffing and MagLev, PH-Symb recovers the correct functional pH structure with close coefficients. In particular, the MagLev Hamiltonian contains the nonseparable product $q^2\phi^2$, showing that we can recover multivariate compositions from univariate edge libraries. For PMSM, PH-Symb identifies a state-dependent interconnection matrix, although the individual pH components are not recovered exactly. Expanding the pH state equation $\dot x=(J-R)\nabla H+Bu$, the true PMSM model gives
\[
\begin{aligned}
\dot{\phi}_d
&=83.34p\phi_q-59.21\phi_d+u_1,\\
\dot{\phi}_q
&=-83.34p\phi_d-14.17p-59.21\phi_q+u_2,\\
\dot p
&=44.74\phi_q-0.217p,
\end{aligned}
\]
whereas PH-Symb gives
\[
\begin{aligned}
\widehat{\dot{\phi}}_d
&=79.84p\phi_q-59.38\phi_d+u_1,\\
\widehat{\dot{\phi}}_q
&=-82.59p\phi_d-14.00p-59.81\phi_q+u_2,\\
\widehat{\dot p}
&=-14.44\phi_d\phi_q+42.12\phi_q-0.606p.
\end{aligned}
\]
The two flux dynamics are accurately recovered, while the momentum dynamics retain an additional coupling term and a less accurate dissipation coefficient.

Figure~\ref{fig:iso-rollouts} compares the state trajectories obtained by simulating the identified models with the ground-truth trajectories. Despite the remaining component-level errors, the simulated states closely follow the reference trajectories over the complete test records.

\begin{figure}[!t]
\centering
\includegraphics[width=\columnwidth]{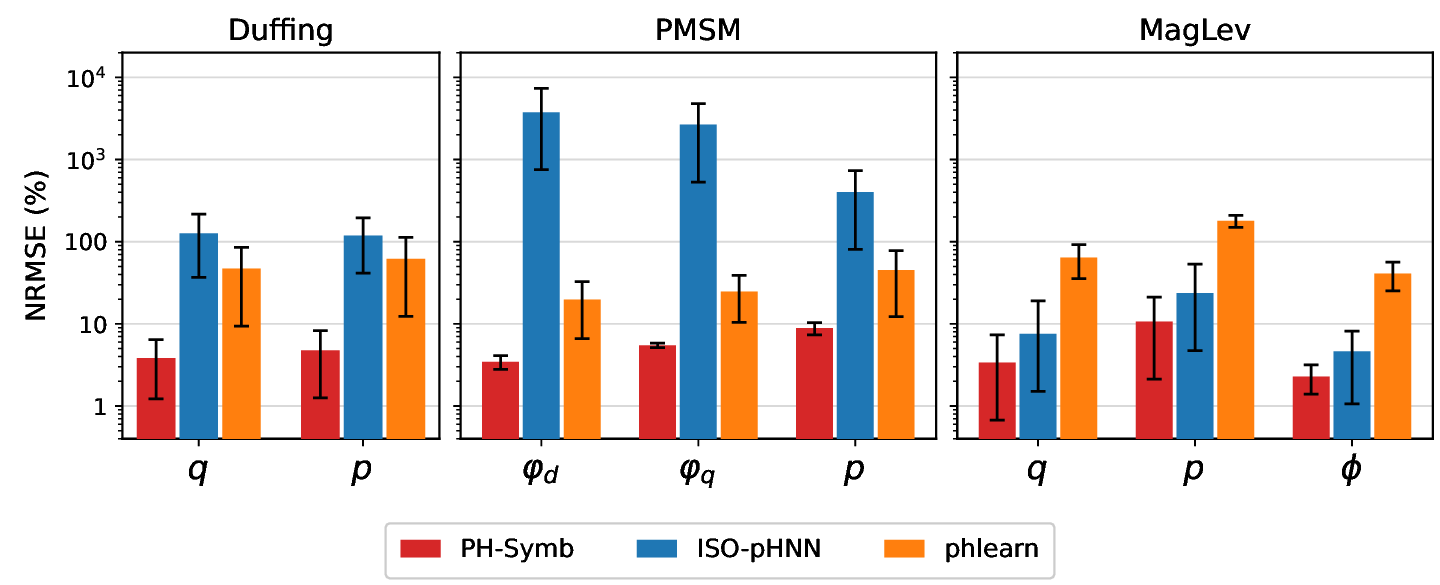}
\caption{State-simulation NRMSE for PH-Symb, ISO-pHNN, and phlearn under matched training-time budgets. Bars and error bars show the mean and standard deviation over ten test experiments. The vertical axis is logarithmic.}
\label{fig:iso-comparison}
\end{figure}

We compare PH-Symb with ISO-pHNN \citep{cherifi2025} and phlearn \citep{eidnes2023} extended with state-dependent MLPs for $J(x)$ and $B(x)$. All methods are trained for approximately $120\,\mathrm{s}$ on Duffing and PMSM and $240\,\mathrm{s}$ on MagLev. Figure~\ref{fig:iso-comparison} reports the mean and standard deviation of the state-simulation NRMSE over ten independent test experiments. PH-Symb obtains the lowest mean NRMSE for every state.

To assess sensitivity to state-measurement noise, independent zero-mean Gaussian noise is added to each training-state channel, with $\sigma$ defined relative to its standard deviation in the clean data. State derivatives are then recomputed using Savitzky--Golay differentiation. Table~\ref{tab:iso-noise} reports the simulation NRMSE after retraining at each noise level.

\begin{table*}[!htbp]
\centering
\caption{True and post-processed PH-Symb models for the synthetic
IO benchmarks. Additive constants in \(H\) are omitted.}
\label{tab:io-models}

\begingroup
\scriptsize
\setlength{\tabcolsep}{0.7pt}
\renewcommand{\arraystretch}{1.22}
\newcommand{\iomatrix}[1]{%
  \left[\begin{smallmatrix}#1\end{smallmatrix}\right]%
}

\begin{tabular*}{\textwidth}{
@{\extracolsep{\fill}}
c cc cc
@{}}
\toprule
&
\multicolumn{2}{c}{\textbf{Nonlinear RLC}}
&
\multicolumn{2}{c}{\textbf{Two mass spring--damper}}
\\
\cmidrule(lr){2-3}
\cmidrule(lr){4-5}
&
\textbf{True}
&
\textbf{PH-Symb}
&
\textbf{True}
&
\textbf{PH-Symb}
\\
\midrule

\(H\)
&
\(\begin{gathered}
\exp(Q_1)-Q_1\\
{}+0.5\phi_1^2+0.5\phi_2^2\\
{}+\exp(Q_2)\\
{}-Q_2
\end{gathered}\)
&
\(\begin{gathered}
0.861\exp(Q_1)-0.862Q_1\\
{}+0.500\phi_1^2+0.503\phi_2^2\\
{}+1.576\exp(Q_2)\\
{}-1.586Q_2-0.110Q_2^2
\end{gathered}\)
&
\(\begin{gathered}
0.025d_1^4+0.500d_1^2+0.500p_1^2\\
{}+0.500d_2^2+0.500p_2^2
\end{gathered}\)
&
\(\begin{gathered}
0.0253d_1^4+0.500d_1^2+0.500p_1^2\\
{}+0.500d_2^2+0.500p_2^2
\end{gathered}\)
\\
\addlinespace[4pt]

\(J\)
&
\(\iomatrix{
0&1&-1&0\\
-1&0&0&0\\
1&0&0&-1\\
0&0&1&0
}\)
&
\(\iomatrix{
0&1.079&-1.080&0\\
-1.079&0&0&0\\
1.080&0&0&-0.862\\
0&0&0.862&0
}\)
&
\(\iomatrix{
0&1&0&0\\
-1&0&-1&0\\
0&1&0&-1\\
0&0&1&0
}\)
&
\(\iomatrix{
0&1.010&0&0\\
-1.010&0&-1.012&0\\
0&1.012&0&-0.998\\
0&0&0.998&0
}\)
\\
\addlinespace[4pt]

\(R\)
&
\(\iomatrix{
1.41\times10^{-5}&0&0&0\\
0&0.707&0&0\\
0&0&0.707&0\\
0&0&0&1.41\times10^{-5}
}\)
&
\(\iomatrix{
0&0&0&0\\
0&0.714&0&0\\
0&0&0.711&0\\
0&0&0&0
}\)
&
\(\iomatrix{
0&0&0&0\\
0&1&0&-0.5\\
0&0&0&0\\
0&-0.5&0&0.5
}\)
&
\(\iomatrix{
0&0&0&0\\
0&0.963&0&-0.465\\
0&0&0&0\\
0&-0.465&0&0.471
}\)
\\
\addlinespace[4pt]

\(B\)
&
\(\iomatrix{
0\\
1\\
0\\
0
}\)
&
\(\iomatrix{
0\\
1.421\\
0\\
0
}\)
&
\(\iomatrix{
0\\
0\\
0\\
1
}\)
&
\(\iomatrix{
0\\
0\\
0\\
0.996
}\)
\\
\bottomrule
\end{tabular*}
\endgroup
\end{table*}

At $30\%$ noise, all mean state NRMSE values remain below $14\%$ for Duffing and $8\%$ for PMSM. MagLev is more sensitive, with errors increasing sharply from $20\%$ noise, especially for $p$ and $\phi$. Inspection of the identified Hamiltonians shows that, at these noise levels, PH-Symb no longer recovers the multiplicative coupling $q^2\phi^2$, which explains the loss of simulation accuracy.

Overall, these results show that PH-Symb can identify explicit symbolic expressions solely from data while preserving the pH structure, even in the presence of nonseparable nonlinearities. However, because the pH decomposition is not unique, accurate dynamics do not necessarily imply recovery of the nonlinearities of the individual pH components.

\begin{table}[!htbp]
\centering
\caption{Free-run state NRMSE (\%) after training with noisy state measurements. Entries are the mean $\pm$ standard deviation over ten test experiments; means are bold.}
\label{tab:iso-noise}
\scriptsize
\setlength{\tabcolsep}{2.5pt}
\renewcommand{\arraystretch}{1.12}
\newcommand{\nrmseentry}[2]{
$\mathbf{#1}\,{\scriptstyle\pm #2}$
}
\begin{tabular}{@{}lcccc@{}}
\toprule
State / $\sigma$
& $0\%$
& $10\%$
& $20\%$
& $30\%$
\\
\midrule
\multicolumn{5}{@{}l}{\emph{Duffing}}
\\
$q$
& \nrmseentry{4.4}{2.1}
& \nrmseentry{5.9}{2.4}
& \nrmseentry{5.1}{2.1}
& \nrmseentry{11.2}{4.6}
\\
$p$
& \nrmseentry{5.3}{2.5}
& \nrmseentry{7.1}{3.0}
& \nrmseentry{6.3}{2.4}
& \nrmseentry{13.6}{6.1}
\\
\addlinespace[2pt]
\multicolumn{5}{@{}l}{\emph{PMSM}}
\\
$\phi_d$
& \nrmseentry{3.8}{0.4}
& \nrmseentry{3.7}{0.4}
& \nrmseentry{3.2}{0.3}
& \nrmseentry{4.2}{0.2}
\\
$\phi_q$
& \nrmseentry{6.0}{0.2}
& \nrmseentry{6.0}{0.2}
& \nrmseentry{6.0}{0.2}
& \nrmseentry{7.1}{0.3}
\\
$p$
& \nrmseentry{9.4}{1.4}
& \nrmseentry{8.8}{1.2}
& \nrmseentry{6.8}{0.9}
& \nrmseentry{7.5}{1.5}
\\
\addlinespace[2pt]
\multicolumn{5}{@{}l}{\emph{MagLev}}
\\
$q$
& \nrmseentry{5.3}{4.8}
& \nrmseentry{5.7}{4.1}
& \nrmseentry{31.8}{34.9}
& \nrmseentry{30.3}{35.2}
\\
$p$
& \nrmseentry{15.3}{11.7}
& \nrmseentry{18.5}{10.5}
& \nrmseentry{74.1}{28.4}
& \nrmseentry{75.4}{27.3}
\\
$\phi$
& \nrmseentry{7.2}{1.9}
& \nrmseentry{4.3}{1.1}
& \nrmseentry{74.7}{8.2}
& \nrmseentry{97.6}{14.2}
\\
\bottomrule
\end{tabular}
\end{table}

\subsection{Input--output benchmarks}
\subsubsection{Synthetic data}

We now evaluate PH-Symb on two systems using simulated IO data. We consider the $N=2$ nonlinear RLC ladder from \citep{chaturantabut2016}, with state $(Q_1,\phi_1,\phi_2,Q_2)$, where $Q_1$ and $Q_2$ are the capacitor charges and $\phi_1$ and $\phi_2$ are the inductor magnetic fluxes, and left-port voltage and current as input and output. We also consider the two-mass spring--damper system (2MSD) system from \citep{moradi2026}, with state $(d_1,p_1,d_2,p_2)$, where $d_1=q_1$ and $d_2=q_1-q_2$. The input acts on the second mass and $y=p_2$. The corresponding ground-truth pH components are reported in Table~\ref{tab:io-models}.

Both systems are sampled at $0.1\,\mathrm{s}$ and excited by the random-phase multisine inputs defined in \eqref{eq:iso-input}, spanning $0.01$--$1\,\mathrm{Hz}$. For each system, the training and test sets contain 20 periods each, whereas the validation set contains 8 periods. Each period contains 1000 samples. Zero-mean Gaussian noise with standard deviation $10^{-3}$ is added to the output. The input RMS is 1.5 for the RLC ladder and 1.58 for the 2 MSD system.

For both benchmarks, we set the state dimension to $n_z=4$, the encoder output- and input-history lengths to $n_a=n_b=20$, and the prediction horizon to $N_p=150$.

\begin{figure*}[!t]
\centering

\begin{minipage}[t]{0.48\textwidth}
  \centering
  \includegraphics[width=\linewidth]{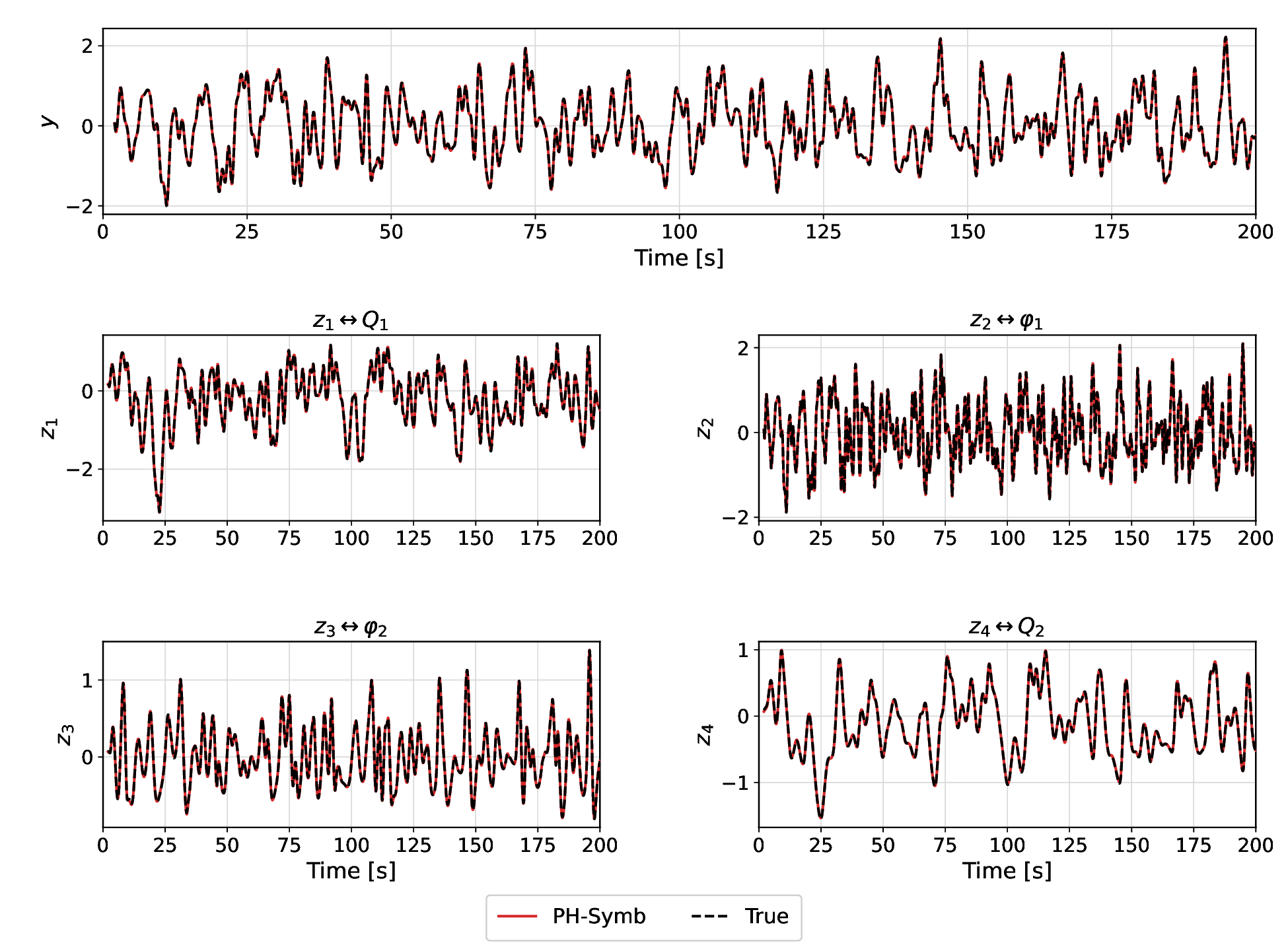}
  \vspace{-1mm}

  \small (a) Nonlinear RLC 
\end{minipage}
\hfill
\begin{minipage}[t]{0.48\textwidth}
  \centering
  \includegraphics[width=\linewidth]{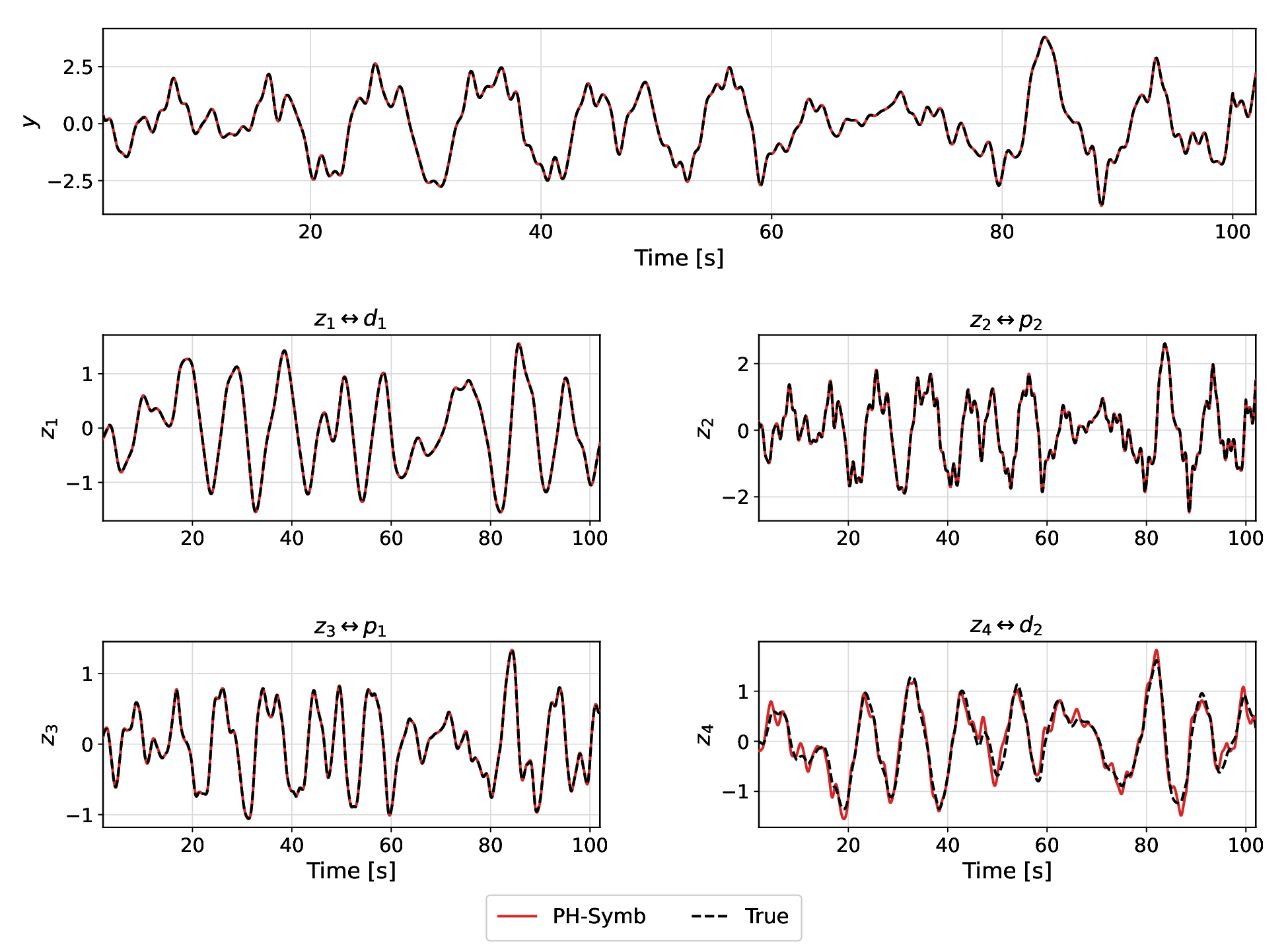}
  \vspace{-1mm}

  \small (b) Two-mass spring--damper system
\end{minipage}

\caption{Simulation results for the synthetic IO benchmarks. The upper panel of each plot shows the predicted output, while the lower panels show the state trajectories generated by the identified pH model.}
\label{fig:io-synthetic-rollouts}
\end{figure*}

The models obtained after training have low output-simulation errors, as shown in Table~\ref{tab:io-models}, but their matrices $J$, $R$, and $B$ are dense. To expose a sparse structure, the internal coordinates with quadratic energy terms are centered and scaled so that each contributes the same quadratic term to $H$. Their total contribution to $H$ is then preserved by any orthogonal mixing of these coordinates. We select the orthogonal transformation that gives the sparsest representation of $J$, $R$, and $B$.

For the RLC ladder, the post-processed model recovers the input acting on a single flux coordinate and assigns dissipation to the two flux coordinates. The capacitor energies have the expected exponential form, although the second retains a quadratic correction. For the 2MSD system, the dissipation acts on the momentum coordinates and the input acts on $p_2$. The quartic term is assigned to the nonlinear spring coordinate.Overall, the post-processed models recover the expected physical structure of both systems, although the quadratic correction shows that the RLC components are not recovered exactly.

Figure~\ref{fig:io-synthetic-rollouts} shows that PH-Symb accurately reproduces both the measured output and the coordinate-aligned internal trajectories, despite only input--output data being available for training. This coordinate alignment is possible here because the state dimension is assumed to be known and the Hamiltonian is separable. The energy can therefore be examined independently for each coordinate, allowing us to identify the coordinates with quadratic energy and rotate only these coordinates to obtain sparse matrices without changing the Hamiltonian. Without these assumptions, accurate output simulation would not necessarily imply recovery of the physical state coordinates.

\subsubsection{Experimental data}

The previous IO experiments use long and simulated records. We now consider the cascaded-tank benchmark to evaluate PH-Symb on a short experimental record and to compare it with IO system identification methods that do not necessarily preserve the pH structure. The system consists of two tanks connected in series. The input is the pump voltage, and the output is the measured water level in the lower tank \citep{schoukens2017benchmarks}.

The data set contains two records of $1024$ samples, acquired with a sampling period of $4\,\mathrm{s}$. The first record is used for training, and the first $512$ samples of the second record are used for validation. The complete second record is used for evaluation. Inputs and outputs are standardized using the training record, and the errors are reported in volts after returning the output to its original scale. We set $n_a=n_b=9$ and $n_z=4$.

The limited amount of data leads to a compromise between simulation accuracy and model sparsity. We therefore consider two variants. The dense variant retains the complete symbolic expansion, whereas the sparse variant removes terms with negligible contributions and re-estimates the coefficients of the retained terms. We use the prior knowledge that square-root nonlinearities are present in this system. To avoid numerical instability near zero, the library includes the regularized signed-square-root function $\rho(s)=\operatorname{sign}(s)\bigl(\sqrt{|s|+0.05}-\sqrt{0.05}\bigr)$. 

The dense model retains $416$ active library coefficients and is therefore not reported term by term. The sparse model retains $28$ active library coefficients, which combine into the following pH components: 
{\scriptsize
\setlength{\jot}{5pt}
\begin{equation}
\label{eq:tanks-interpretable}
\begin{gathered}
\widehat H(z)
=
0.0283z_1^3+0.0403z_1^2
+0.6800z_2^2+0.5504z_4^2,
\\[2pt]
\widehat J(z)
=
\left[
\begin{smallmatrix}
0 & j_{12} & 0 & j_{14}\\
-j_{12} & 0 & 0 & j_{24}\\
0 & 0 & 0 & 0\\
-j_{14} & -j_{24} & 0 & 0
\end{smallmatrix}
\right],
\qquad
\widehat R(z)
=
\left[
\begin{smallmatrix}
r_1 & 0 & 0 & 0\\
0 & \varepsilon & 0 & 0\\
0 & 0 & \varepsilon & 0\\
0 & 0 & 0 & r_4
\end{smallmatrix}
\right],
\\[2pt]
\begin{aligned}
j_{12}
&=
0.3412z_2
+0.3586\rho(z_2)
+0.1940\rho(z_4),
\\
j_{14}
&=
0.0014z_3^2
-0.1301\rho(z_4),
\\
j_{24}
&=
-0.0595z_3^2
-0.3123,
\\
r_1
&=
\left(0.0498z_3^2-0.3271\right)^2,
\\
r_4
&=
\left(0.2507z_4^2-0.0826\right)^2,
\\
\widehat B(z)
&=
\begin{bmatrix}
b_1 & b_2 & 0 & b_4
\end{bmatrix}^{\top},
\\
b_1
&=
0.0677z_2^2
-0.2301z_2
-0.4065\rho(z_2)
-0.1811\rho(z_4),
\\
b_2
&=
0.0568z_4^2
-0.0403\rho(z_2),
\\
b_4
&=
0.0419z_3^2
+0.1211z_4
+0.4384.
\end{aligned}
\end{gathered}
\end{equation}
}

These expressions describe the identified pH dynamics in learned coordinates, so their coefficients do not correspond directly to the physical tank parameters.

\begin{figure}[!htbp]
\centering
\includegraphics[width=\columnwidth]{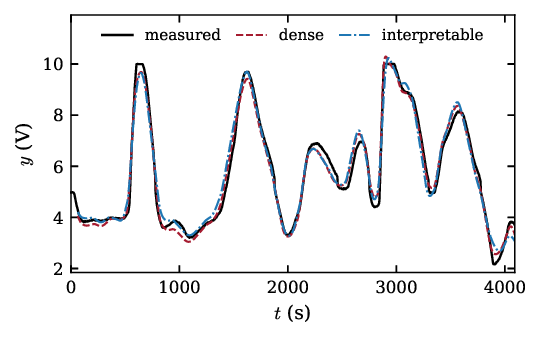}
\caption{Measured cascaded-tank output and free-run predictions of the dense and interpretable PH-Symb variants on the evaluation record of the benchmark}
\label{fig:tanks-output}
\end{figure}

Table~\ref{tab:tanks-comparison} compares the dense and sparse PH-Symb variants with published IO identification methods that either impose a pH structure or use unconstrained model classes. The dense variant achieves an output-simulation RMS error of $0.26\,\mathrm{V}$ and remains close to the best reported results. The sparse variant is less accurate, with an RMS error of $0.35\,\mathrm{V}$, but provides a much sparser symbolic model and still follows the measured output closely, as shown in Figure~\ref{fig:tanks-output}. Overall, the same PH-Symb framework can be used to favor simulation accuracy by retaining the complete symbolic expansion, or favor sparsity by removing terms with small contributions.

\begin{table}[!t]
\centering
\caption{Cascaded-tank free-run output RMS; lower is better.}
\label{tab:tanks-comparison}
\begingroup
\footnotesize
\setlength{\tabcolsep}{3pt}
\renewcommand{\arraystretch}{1.05}
\begin{tabular*}{\columnwidth}{
@{\extracolsep{\fill}}
lc
@{}}
\toprule
Method & RMS (V) \\
\midrule
CT subspace encoder \citep{beintema2023}
& \(0.22\) \\
\textbf{Dense PH-Symb}
& \(\mathbf{0.26}\) \\
OE-pHNN \citep{moradi2026}
& \(0.28\) \\
Tensor B-splines \citep{karagoz2020tensor}
& \(0.30\) \\
TSEM \citep{forgione2021continuous}
& \(0.33\) \\
\textbf{Sparse PH-Symb}
& \(\mathbf{0.35}\) \\
DT subspace encoder \citep{beintema2023}
& \(0.37\) \\
SCI \citep{forgione2021continuous}
& \(0.40\) \\
State-space GP prior \citep{svensson2017}
& \(0.45\) \\
LSTM \citep{champneys2024}
& \(0.49\) \\
LTI state-space \citep{champneys2024}
& \(0.59\) \\
\bottomrule
\end{tabular*}
\endgroup
\end{table}

\section{Conclusion}
\label{sec:conclusion}
This paper introduced PH-Symb, a framework for identifying symbolic port-Hamiltonian models from both ISO and IO data. By combining pH-preserving KANs with edge-wise sparse regression, PH-Symb obtains explicit expressions for the pH components while preserving the interconnection, dissipation, and port structure. In addition, a new normalization was proposed to preserve the pH structure  and  thus express the identified model in the original coordinates. Synthetic and experimental results show that PH-Symb provides interpretable models with lower errors than the neural baselines. Future work will address high-dimensional and distributed-parameter systems and extend the method to irreversible pH models.

\bibliographystyle{elsarticle-harv}       
\bibliography{autosam} 

\appendix

\end{document}